\documentclass[a4paper,12pt]{article}
\usepackage{latexsym}
\usepackage{amsfonts}
\newtheorem{theorem}{\sc Theorem}[section]
\newtheorem{corollary}[theorem]{\sc Corollary}

\newtheorem{remark}[theorem]{\sc Remark}
\newtheorem{proposition}[theorem]{\sc Proposition}

\title{\sc A remark on $\beta$-locally closed sets}
\author{Julian {\sc Dontchev} and Maximilian {\sc Ganster}}
\date{}
\begin{document}
\baselineskip=20pt plus 1pt minus 1pt
\newcommand{\bl}{$\beta$-locally}
\newcommand{\bs}{$\beta$-submaximal}
\newcommand{\fxy}{$f \colon (X,\tau) \rightarrow (Y,\sigma)$}
\maketitle
\begin{abstract}
The aim of this note is to show that every subset of a given
topological space is the intersection of a preopen and a
preclosed set, therefore \bl \ closed, and that every topological
space is \bs.
\end{abstract}

\section{Introduction}\label{s1}

In a recent paper, Gnanambal and Balachandran  \cite{GB1}
introduced the classes of \bl\ closed sets, \bs\ spaces and
$\beta$-LC-continuous functions. The purpose of our note is to
show that every subset of any topological space is the
intersection of a preopen set and a preclosed set, hence \bl\
closed, and therefore every function \fxy\ is
$\beta$-LC-continuous.  We have felt the need to point out
explicitly this observation since over the years several papers
have investigated concepts like "pre-locally closed sets" or
"\bl\ closed sets" which do not have any nontrivial meaning. In
addition, we will show that every space is \bs\ and we will point
out that most results of \cite{GB1} are either trivial or false.

Let $A$ be a subset of a topological space $(X,\tau)$.
Following Kronheimer \cite{K1}, we call the interior of the
closure of $A$, denoted by $A^+$, the {\em consolidation} of $A$.
Sets included in their consolidation are called {\em preopen}
or {\em locally dense}. Complements of preopen sets are called
{\em preclosed} and the preclosure of a set $A$, denoted by ${\rm
pcl} (A)$, is the intersection of all preclosed supersets of $A$.
Since union of preopen sets is also preopen, the preclosure
of every set is in fact a preclosed set. If $A$ is included
in the closure of its consolidation, then $A$ is called {\em
$\beta$-open} or {\em semi-preopen}. Complements of $\beta$-open
sets are called {\em $\beta$-closed}. The $\beta$-closure of $A$,
denoted by ${\rm cl}_{\beta} (A)$ is the intersection of all
$\beta$-closed supersets of $A$. In \cite{GB1}, Gnanambal and
Balachandran called a set $A$ {\em \bl\ closed} if $A$ is
intersection of a $\beta$-open and a $\beta$-closed set. They
defined a set $A$ to be {\em $\beta$-dense} \cite{GB1} if ${\rm
cl}_{\beta} (A) = X$ and called a space $X$ {\em \bs\ }
\cite{GB1} if every $\beta$-dense subset is $\beta$-open. A
function \fxy\ is called {\em $\beta$-LC-continuous} \cite{GB1}
if the preimage of every open subset of $Y$ is \bl\ closed in
$X$.

The following implications hold and none of them is reversible:

\begin{center}
dense $\Rightarrow$ preopen $\Rightarrow$ $\beta$-open
$\Rightarrow$ \bl\ closed
\end{center}

\section{Every set is \bl\ closed}\label{s2}

\begin{proposition}\label{p1}
Every subset $A$ of a topological space $(X,\tau)$ is the
intersection of a preopen and a preclosed set, hence pre-locally
closed.
\end{proposition}

{\em Proof.} Let $A \subseteq (X,\tau)$. Set $A_1 = A \cup (X
\setminus {\rm cl} (A))$. Since $A_1$ is dense in $X$, it is also
preopen. Let $A_2$ be the preclosure of $A$, i.e., $A_2 = A \cup
{\rm cl} ({\rm int} (A))$. Clearly, $A_2$ is a preclosed set.
Note now that $A = A_1 \cap A_2$. $\Box$

\begin{corollary}\label{c1}
{\rm (i)} Every set is \bl\ closed and every function is
$\beta$-LC-continuous.

{\rm (i)} Every topological space is \bs.
\end{corollary}

{\em Proof.} (i) Every preopen (resp.\ preclosed) set is
$\beta$-open (resp.\ $\beta$-closed).

(ii) By \cite[Corollary 3.24]{GB1} a topological space is \bs\
if and only if every set is \bl\ closed.

\begin{remark}\label{r1}
{\em (i) Corollary~\ref{c1} makes \cite{GB1} trivial.

(ii) Example 3.4 from \cite{GB1} is wrong as the subset $A = \{
\frac{1}{n} \colon n = 1,2,\ldots \} \cup (2,3) \cup (3,4) \cup
\{ 4 \} \cup (5,6) \cup \{ x \colon x$ is irrational and $7 \leq
x < 8 \}$ of the real line $\mathbb R$ is indeed \bl\ closed.

(iii) Proposition 3.6 from \cite{GB1} is wrong as every proper
nonempty subset of the real line $\mathbb R$ with the indiscrete
topology is $\beta$-open and preclosed but not semi-open.

(iv) Example 4.11 from \cite{GB1} is wrong, since the space
$(X,\tau)$, where $X = \{ a,b,c,d \}$, $\tau = \{ \emptyset, \{
a,b \}, \{ c,d \}, X \}$ is {\em not} an ${\alpha}{\beta}$-space.
Note that $\{ a \}$ is $\beta$-open but not $\alpha$-open (an
{\em $\alpha$-open set} is a set which is the difference of an
open and a nowhere dense set).

(v) An ${\alpha}{\beta}$-space \cite{GB1} is in fact a strongly
irresolvable, extremally disconnected space.

(vi) An $\alpha$-locally closed set (\cite[Definition 2.1
(x)]{GB1}) is nothing else but a simly-open set.}
\end{remark}

\
\begin{center}
Department of Mathematics\\University of Helsinki\\PL 4,
Yliopistonkatu 5\\00014 Helsinki\\Finland\\e-mail: {\tt
dontchev@cc.helsinki.fi}
\end{center}
\
\begin{center}
Department of Mathematics\\Graz University of
Technology\\Steyrergasse 30\\A-8010 Graz\\Austria\\e-mail:
{\tt ganster@weyl.math.tu-graz.ac.at}
\end{center}
\
\end{document}